\documentclass[12pt]{article}
\usepackage{mathrsfs}
\usepackage{stmaryrd}
\usepackage{amsfonts,amsmath,amssymb,amscd}
\usepackage{shadow}
\usepackage{graphicx}
\usepackage{color}
\usepackage{longtable}
\usepackage{pstricks,multido}
\usepackage{hyperref}
\usepackage{qtree}
\allowdisplaybreaks

\newtheorem{thm}{Theorem}[section]
\newtheorem{lem}[thm]{Lemma}

\newtheorem{conj}{Conjecture}[section]
\def\qed{\hfill \rule{4pt}{7pt}}

\def\pf{\noindent {\it{Proof.} \hskip 2pt}}

\numberwithin{equation}{section}
\numberwithin{figure}{section}

\begin{document}
\begin{center}
{\large\bf On Pattern Avoiding Alternating Permutations}
\end{center}

\begin{center}
Joanna N. Chen$^1$, William Y.C. Chen$^2$, Robin D.P. Zhou$^3$

$^{1,3}$Center for Combinatorics, LPMC-TJKLC\\
Nankai University\\
Tianjin 300071, P.R. China

$^2$Center for Applied Mathematics\\
Tianjin University\\
Tianjin 300072, P.R. China

$^1$joanna@cfc.nankai.edu.cn, $^2$chenyc@tju.edu.cn, $^3$robin@cfc.nankai.edu.cn % $^{\wedge}$ .

\end{center}

\begin{abstract}
An alternating permutation of length $n$ is a permutation $\pi=\pi_1 \pi_2 \cdots \pi_n$ such that $\pi_1 < \pi_2 > \pi_3 < \pi_4 > \cdots$.
 Let $A_n$ denote set of alternating permutations of $\{1,2,\ldots, n\}$, and let $A_n(\sigma)$ be set of alternating permutations
 in $A_n$ that avoid a pattern $\sigma$.
Recently, Lewis used  generating trees to enumerate $A_{2n}(1234)$, $A_{2n}(2143)$ and $A_{2n+1}(2143)$, and he posed several
conjectures on the Wilf-equivalence of alternating permutations
avoiding certain patterns. Some of these conjectures have been proved by B\'ona, Xu and Yan.
In this paper, we  prove the two relations  $|A_{2n+1}(1243)|=|A_{2n+1}(2143)|$ and $|A_{2n}(4312)|=|A_{2n}(1234)|$ as conjectured by Lewis.
\end{abstract}



\section{Introduction}

 The objective of this paper is to prove two conjectures of
 Lewis on the Wilf-equivalence of alternating permutations
 avoiding certain patterns of length four.

  We begin with some notation and terminology.
  Let $[n]=\{1,2,\cdots, n\}$, and let $S_n$
  be the set of permutations of $[n]$. A permutation $\pi=\pi_1 \pi_2 \cdots \pi_n$ is said to be an alternating permutation
  if $\pi_1 < \pi_2 > \pi_3 < \pi_4 > \cdots$. An alternating permutation is also called
  an up-down permutation. A permutation  $\pi$ is said to be a down-up
  permutation if $\pi_1 > \pi_2 < \pi_3 > \pi_4 < \cdots$. We denote by $A_n$ and $A'_n$  the set of  alternating permutations
  and the set of down-up permutations  of $[n]$, respectively.
  For a permutation
   $\pi \in S_n$, its reverse $\pi^r \in S_n$ is defined by
  $\pi^r(i)=\pi(n+1-i)$ for $1\leq i\leq n$. The complement  of $\pi$, denoted $\pi^c \in S_n$,  is defined by $\pi^c(i)=n+1-\pi(i)$ for  $1\leq i\leq n$.
  It is clear that the complement  operation
   gives a bijection between $A_n$ and $A'_n$.

  Given a permutation $\pi$ in $S_n$ and  a permutation $\sigma=\sigma_1 \sigma_2 \cdots \sigma_k \in S_k$, where $k\leq n$, we say
 that $\pi$ contains a pattern $\sigma$  if there exists a subsequence
 $\pi_{i_1}\pi_{i_2} \cdots \pi_{i_k}$ $(1 \leq i_1 < i_2 < \cdots <i_k \leq n)$ of $\pi$ that is order isomorphic to $\sigma$,
 in other words,  for all $l,m \in [k]$, we have
 $\pi_{i_l}<\pi_{i_m}$ if and only if $\sigma_l <\sigma_m$. Otherwise, we say that $\pi$  avoids a pattern $\sigma$, or
 $\pi$ is $\sigma$-avoiding. For example, $74538126$ is 1234-avoiding, while it
  contains pattern $3142$ corresponding to the subsequence $7486$.

  Let $S_n(\sigma)$ denote the set of permutations of length $n$ that  avoid a pattern $\sigma$. Let $A_n(\sigma)$
 denote the set of $\sigma$-avoiding alternating permutations of $[n]$, and let $A'_n(\sigma)$ denote the set of
 $\sigma$-avoiding down-up permutations of $[n]$.
 Mansour \cite{Mansor2003} showed that $|A_{2n}(132)| = C_n$, where
 $C_n$ is the Catalan number
\[ \frac{1}{n+1}{2n \choose n}.\]
Meanwhile, Deutsch and Reifegerste ( as reported by Stanley \cite{rp2011}) showed  that  $|A_{2n}(123)| = C_n$.
Recently, Lewis \cite{Lewis2012} showed that the generating trees for $A_{2n}(1234)$ and $A_{2n}(2143)$  are isomorphic  to the
generating tree for the set of standard Young tableaux of shape $(n,n,n)$. From the hook-length formula it follows that
\begin{equation} \label{V_n}
 |A_{2n}(1234)|=|A_{2n}(2143)|=\frac{2(3n)!}{n!(n+1)!(n+2)!}.
 \end{equation}
The above number is called the $n$-th $3$-dimensional Catalan number, and we
shall denote it by $C_n^{(3)}$. Notice that $C_n^{(3)}$ also equals the number of
walks in $3$-dimensions from $(0,0,0)$ to $(n,n,n)$ by using steps $(1,0,0)$, $(0,1,0)$ and $(0,0,1)$ that do not go below the plane $x=y=z$.
  Lewis  showed that  $A_{2n+1}(2143)$ has the same generating
  tree as that of shifted standard Young tableaux of shape $(n+2,n+1,n)$. Using the hook-length formula for shifted standard Young tableaux given by Krattenthaler \cite{C.Krattenthaler}, we deduce that
 \[ |A_{2n+1}(2143)|=\frac{2(3n+3)!}{n!(n+1)!(n+2)!(2n+1)(2n+2)(2n+3)}.\]

  The following conjectures were posed by Lewis \cite{Lewis2012}.

  \begin{conj}
  For $n\geq 1$ and  $\sigma \in \{ 1243, 2134, 1432, 3214, 2341,4123, 3421, 4312\}$, we have
  \[|A_{2n}(\sigma)|=|A_{2n}(1234)|=|A_{2n}(2143)|.\]
  \end{conj}

  \begin{conj}
  For $n\geq 0$ and  $\sigma \in \{ 2134,$   $4312, $ $ 3214, $ $4123\}$, we have
  \[|A_{2n+1}(\sigma)|=|A_{2n+1}(1234)|.\]
  \end{conj}

  \begin{conj}
   For $n\geq 0$ and  $\sigma \in \{ 1243, $  $3421, $  $1432, $  $ 2341\}$, we have
  \[|A_{2n+1}(\sigma)|=|A_{2n+1}(2143)|.\]
   \end{conj}

By showing that a classical bijection on pattern
 avoiding permutations preserves the alternating property,
 B\'ona \cite{BoNA} proved that
\begin{align}
   |A_{2n}(1243)|&=|A_{2n}(1234)|,\label{c3}\\
    |A_{2n+1}(2134)|&=|A_{2n+1}(1234)|.
\end{align}
   Xu and Yan \cite{Yan2012}
  constructed bijections that lead to the following relations
  \begin{align*}
   |A_{2n}(4123)|=|A_{2n} &(1432)|=|A_{2n}(1234)|,\\
   |A_{2n+1}(1432)|&=|A_{2n+1}(2143)|,\\
   |A_{2n+1}(4123)|&=|A_{2n+1}(1234)|.
  \end{align*}
As for the above conjectures, there are essentially two
unsolved cases, namely, \begin{equation}\label{c1}
|A_{2n+1}(1243)|=|A_{2n+1}(2143)|,\end{equation}
and
\begin{equation}\label{c2}
|A_{2n}(4312)|=|A_{2n}(1234)|,
\end{equation}
 because the other
remaining cases can be deduced by the reverse and complement  operations.

In this paper, we prove the above conjectures (\ref{c1}) and (\ref{c2}). To be more
specific, we show that the generating tree for $A_{2n+1}(1243)$
 coincides with the generating tree for $A_{2n+1}(2143)$ as given
 by Lewis \cite{Lewis2012}. So we are led to relation (\ref{c1}). The construction of the generating tree
 for $A_{2n+1}(1243)$ can be adapted to obtain the
generating tree for $A_{2n}(1243)$, which turns out to be isomorphic to the generating tree for $A_{2n}(1234)$ as constructed by Lewis \cite{Lewis2012}. This gives another proof of relation (\ref{c3})
conjectured by Lewis and proved by B\'ona.

To prove (\ref{c2}), we show that
the generating tree for $A_{2n+1}(1243)$ is isomorphic to the generating tree for the set of shifted standard Young tableaux of shape $(n+2,n+1,n)$ as given by Lewis \cite{Lewis2012}. We adopt the notation  $SHSYT(\lambda)$ for the set of shifted standard Young tableaux of shape $\lambda$.
As can be easily seen, a label $(a,b)$ in the generating tree for $A_{2n+1}(1243)$ corresponds to a label $(a+1,b)$  in the generating tree for $SHSYT(n+2,n+1,n)$. By restricting the correspondence to   certain labels of the generating trees,
we obtain a bijection between a subset  of $A_{2n+1}(1243)$ and a subset of $SHSYT(n+2,n+1,n)$. This leads to the relation $|A_{2n}(4312)|=|SHSYT(n+2,n,n-2)|$. By the hook-length formula for shifted standard Young tableaux, we see that $SHSYT(n+2,n,n-2)$ is counted by $C_n^{(3)}$. Since $A_{2n}(1234)$ is also enumerated by $C_n^{(3)}$,
 we arrive at relation (\ref{c2}).

Since we already have that $|A_{2n}(4312)|=|A_{2n}(1234)|$, it is natural to consider whether one can construct a generating tree for $A_{2n}(4312)$ that is
isomorphic to the generating tree for $A_{2n}(1234)$ given by Lewis.
While we have not found such a generating tree for $A_{2n}(4312)$,
we obtain a generating tree for $A_{2n}(4312)$ that
can be used to give a second proof of
 relation (\ref{c1}).
By deleting the leaves of the generating tree for
$A_{2n}(4312)$ and changing the label $(a,b)$ to $(a-1,b)$,
 we are led to the generating  tree for $A_{2n}(3412)$ as given by Lewis \cite{Lewis2012}. Furthermore, by restricting
 this correspondence to certain labels, we obtain
 relation (\ref{c1}).

This paper is organized as follows.
In Section \ref{sec:1243}, we construct a generating tree for
$A_{2n+1}(1243)$, which turns out to be the same with the generating tree for $A_{2n+1}(2143)$ given by Lewis. This proves
relation (\ref{c1}). By similar
constructions, we see that $A_{2n}(1234)$ and $A_{2n}(1243)$
have isomorphic generating trees. This yields another proof of (\ref{c3}). In Section \ref{sec:4312}, we  prove
(\ref{c2}) by showing that $|A_{2n}(4312)|$ is equal to the number of shifted standard Young tableaux of shape $(n+2,n,n-2)$.
In Section \ref{sec:4312tree}, we construct a generating tree for $A_{2n}(4312)$, and give another proof of (\ref{c1}).

\section{Generating trees for $A_{2n+1}(1243)$ and $A_{2n}(1243)$} \label{sec:1243}

In this section, we construct the generating tree for
$A_{2n+1}(1243)$ which turns out to be the same as
that for $A_{2n+1}(2143)$. This proves (\ref{c1}), which we restate as the
following theorem.

\begin{thm}\label{thm:1243=2143}
For $n\geq 0$, we have
$|A_{2n+1}(1243)|=|A_{2n+1}(2143)|$.
\end{thm}

We also obtain a generating tree for $A_{2n}(1243)$
and we show that it is isomorphic to the generating tree for $A_{2n}(1234)$.
This confirms (\ref{c3}), which we restate as the following theorem.

\begin{thm}\label{thm:1243=1234}
For $n\geq 1$, we have
$|A_{2n}(1243)|=|A_{2n}(1234)|$.
\end{thm}

Let us give an overview of the terminology on generating trees. Given a sequence $\{\Sigma_n\}_{n \geq 1}$ of finite, nonempty sets with $|\Sigma_1|=1$, a generating tree for
this sequence is a rooted, labeled tree such that the vertices at level $n$ are the elements of $\Sigma_n$ and the label of each vertex determines the multiset of labels of its children.
Thus, the generating tree is fully described by its root vertex and the succession rule $L \to S$ which gives the set $S$ of labels of
the children in terms of the label $L$ of their parent. Here, we denote a generating tree in the following form,
\begin{align*}
    \left\{
      \begin{array}{ll}
        root \colon & \hbox{\text{the label of the root}}, \\ [3pt]
        rule \colon & \hbox{\text{succession rules}}.
      \end{array}
    \right.
\end{align*}
Sometimes we also refer a generating tree for $\Sigma_n$ to the
generating tree for the sequence $\{\Sigma_n\}_{n\geq 1}$.

The construction of a generating tree for $\{\Sigma_n\}_{n\geq1}$
requires the generation of $\Sigma_{n+1}$ based
on $\Sigma_{n}$. For $u\in \Sigma_n$, $w\in \Sigma_{n+1}$, let $w$ be a child of $u$ in the generating tree if and only if
$w$ is generated by $u$. Thus it is sufficient to determine the structure of the generating tree by defining the children of each element.

To illustrate the idea of generating trees,
we consider the construction of a generating tree for $S_n$.
We need to determine the children of each permutation in $S_n$.
Given $\pi \in S_{n}$, we can generate $n+1$ permutations
in $S_{n+1}$. For $ 1\leq i  \leq n+1$,
let $i \mapsto \pi$ denote the permutation $\sigma=\sigma_1
\sigma_2\cdots \sigma_{n+1}$  in $S_{n+1}$ such that
$\sigma_1=i$ and $\sigma_2 \sigma_3 \cdots\sigma_{n+1}$ is order
isomorphic to $\pi$. In other words, $i \mapsto \pi$
 is the permutation obtained from $\pi$ by adding
  $i$ to the beginning of $\pi$ and
increasing each element not less than $i$ by $1$.
For example, $3 \mapsto 3142 = 34152$
is a child of $3142$ in the generating tree.

Notice that Lewis \cite{Lewis2012} used the
 notation $\pi \leftarrow i$ denote the permutation $\sigma=\sigma_1
\sigma_2\cdots \sigma_{n+1}$  in $S_{n+1}$ such that
 $\sigma_{n+1}=i$ and
$\sigma_1 \sigma_2 \cdots \sigma_n$ order isomorphic to $\pi$.
The idea of generating trees is to give succession rules for the structure
of the generating tree by assigning labels to the vertices.
For the case of permutations, given $\pi=\pi_1\pi_2 \cdots \pi_n \in S_n$,
we associate it with a label $(\pi_1,n)$. Then we have the
generating tree for $S_n$ as follows
\begin{align*}
    \left\{
      \begin{array}{ll}
        root\colon & \hbox{$(1,1)$}, \\[3pt]
        rule\colon & \hbox{$(i,n)\to \{(j,n+1)\,|\,1\leq j \leq n+1\}$}.
      \end{array}
    \right.
\end{align*}

By the recursive construction of alternating permutations,
Lewis \cite{Lewis2012} obtained generating schemes for $A_{2n}$
and $A_{2n}(\sigma)$. Here we describe the recursive
constructions of $A_{2n}$ and  $A_{2n}(\sigma)$
by adding elements at the beginning. This choice of notation
seems to be more convenient for the description of the
construction of the generating trees for $A_{2n+1}(1243)$ and $A_{2n}(4312)$.

For $n\geq 1$, let $u=u_1 u_2 \cdots u_{2n}$ be an alternating
permutation in $A_{2n}$. The generating tree is constructed based on
the following generating scheme.  Consider alternating permutations $w=w_1w_2w_{3}\cdots w_{2n+2}$  in $A_{2n+2}$ such that $w_3 w_4 \cdots w_{2n+2}$ is order isomorphic to $u$. Such permutations are set to be the children of $u$ in the generating
tree. One can also use this recursive procedure to construct
pattern avoiding alternating permutations.
To be specific, given $u\in A_{2n}(\sigma)$,
the set of the children of $u$ is precisely the set
$\{w~|~w = v_1 \mapsto (v_2 \mapsto u),~ w \in A_{2n+2}(\sigma)\}$.
The generating scheme for
pattern avoiding alternating permutations of odd length can be
constructed in the same manner.

We now  proceed to construct the
generating trees for $A_{2n+1}(1243)$ and $A_{2n}(1243)$.
In fact, these two sets have the same succession rules
with different roots.
Here we shall only present the derivation of the succession rules for $A_{2n+1}(1243)$. To this end, we need
to characterize the set of $1243$-avoiding
alternating permutations in $A_{2n+3}$ that are generated by an
alternating permutation $u$ in $A_{2n+1}(1243)$. Such a characterization
leads to a labeling along with succession rules.
% where a label of an alternating permutation can be viewed as a
%classification with respect to the recursive generation
%of $1243$-avoiding alternating permutations.

\begin{thm}\label{lem:characterize}
For $n \geq 0$, given a permutation $u=u_1u_2\cdots u_{2n+1} \in A_{2n+1}(1243)$, define
\begin{align*}
f(u)&= \text{max} \{0,\,u_j\,|\,\text{there exists}~ i ~\text{such that}~i< j  ~\text{and} ~ u_i > u_j\},\\[6pt]
e(u)&= \text{max} \{0,\,u_i\,|\,\text{there exist}~ j ~\text{and}~ k~ \text{such that}~ i<j<k~ \text{and} ~u_i<u_k<u_j\}.
\end{align*}
Then $w$ is a child of $u$ if and only if
it is of the form $w=v_1 \mapsto (v_2 \mapsto u)$, where
\begin{equation} \label{rr1}
e(u) < v_1 \leq v_2,
\end{equation}
and
\begin{equation}\label{rr2}
\text{max}\{u_1+1,
f(u)+1\} \leq  v_2 \leq 2n+2.
\end{equation}
\end{thm}

\pf
Suppose $w=w_1w_2\cdots w_{2n+3}$ is a child of $u$, by definition,
$w$ is of the form $v_1 \mapsto (v_2 \mapsto u)$ and $w \in A_{2n+3}(1243)$.
 Since $w$ is alternating on $[2n+3]$, we have $v_1 \leq v_2 \leq 2n+2$ and $v_2 \geq u_1+1$. By the order of the insertions of
  $v_1$ and $v_2$, we see that $w_1=v_1$ and $w_2=v_2+1$.
  Since $w$ is 1243-avoiding, we have $v_2 \geq f(u)+1$;
   Otherwise, there exists $i < j$ such that $u_i >u_j$ and
 $v_2 \leq u_j$. This implies that $w_1 w_2 w_{i+2} w_{j+2} = v_1(v_2+1)(u_i+2)(u_j+2)$ forms a $1243$-pattern. Moreover, we have $v_1 >e(u)$; Otherwise, there exist $i< j< k$ such that $u_i < u_k <u_j$ and $v_1 \leq u_i$. Clearly,
 $w_{i+2} > u_i \geq v_1 $.
 Thus $w_1 w_{i+2} w_{j+2} w_{k+2}$ is of pattern $1243$.  So we are
 led to the relations (\ref{rr1}) and
 (\ref{rr2}).

Conversely, we assume that $v_1$ and $v_2$ are integers
satisfying conditions (\ref{rr1}) and (\ref{rr2}). We wish to show that $w= v_1 \mapsto (v_2 \mapsto u)$
is an alternating permutation in $A_{2n+3}(1243)$, that is, $w$ is a
child of $u$ in the generating tree.
It is evident from (\ref{rr1}) and (\ref{rr2}) that
$v_1 \leq v_2$ and $v_2 >u_1$. So we
have $
w_1 < w_2 >w_3$, from which we see that $w$ is alternating
 since $w_3 w_4 \cdots w_{2n+3}$ is order isomorphic to
$u$.

It remains to show that $w$ is  $1243$-avoiding.
Assume to the contrary that $w$ contains a $1243$-pattern, that is,
there exist $t <i < j < k$ such that $w_t w_iw_jw_k$
is of pattern $1243$. We claim that $t= 1 ~\text{or} ~2$.
Otherwise, we assume that $t\geq 3$.
Since $u_{t-2} u_{i-2}u_{j-2}u_{k-2}$
is isomorphic to  $w_t w_iw_jw_k$, we find
$u_{t-2} u_{i-2}u_{j-2}u_{k-2}$ forms a $1243$-pattern,
which is a contradiction. It follows that $t\leq 2$.
If $w_2w_iw_jw_k$ forms a $1243$-pattern, then $w_1w_iw_jw_k$ is also a $1243$-pattern. Hence we can always choose $t=1$.
To prove $w$ is $1243$-avoiding, it is sufficient to show that
 it is impossible for $w_1 w_i w_j w_k$ to be a $1243$-pattern.

%Clearly, $i \geq 2$.
We now assume that $w_1w_iw_j w_k$ is a $1243$-pattern.
If $i=2$, we have $w_2<w_k$. Since $w_2=v_2+1$ and $w_k\leq u_{k-2}+2$,  we get
$v_2 \leq u_{k-2}$.
Note that $u_{j-2}u_{k-2}$ is order isomorphic to $w_j w_k$,  so we have $u_{j-2} >  u_{k-2}$. By the definition of $f(u)$, we find
$u_{k-2}\leq f(u)$. It follows that $v_2 \leq f(u)$,
which contradicts to the fact that $v_2 \geq f(u)+1$.
Hence we have $i >2$.

We now claim that $w_1 \leq u_{i-2}$. Otherwise, we assume $w_1 > u_{i-2}$. Since $w_1=v_1$ and $v_1 \leq v_2$,  we find that $u_{i-2}< v_1 \leq v_2$. By the construction of $w$,
 we have $w_i=u_{i-2}$. This yields $w_i<w_1$, which contradicts to
 the  assumption that
$w_1w_iw_j w_k$ is a $1243$-pattern. This proves the claim.
 Clearly, $u_{i-2} u_{j-2} u_{k-2}$ is a $132$-pattern since it is
 order isomorphic to $w_iw_j w_k$. By the definition of $e(u)$, we get
$ u_{i-2}\leq e(u)$. Thus $v_1=w_1  \leq u_{i-2}\leq e(u)$, which contradicts to the fact $v_1 > e(u)$. So we reach the conclusion that
the assumption that $w_1w_iw_j w_k$ is a $1243$-pattern is not valid. In other words, $w$ is $1243$-avoiding. This completes the proof.\qed

In light of the above characterization, we are led to a labeling scheme for
alternating permutations in $A_{2n+1}(1243)$.
For $u \in A_{2n+1}(1243)$, we associate  a label $(a,b)$ to $u$, where
\begin{align}
    &a=2n+2-\text{max}\{u_1+1,f(u)+1\}, \label{labelscheme1}\\[3pt]
    &b=2n+2-e(u). \label{labelscheme2}
\end{align}
For example, the permutation $1 \in A_{1}(1243)$ has label $(0,2)$,
and the  permutation $2546173 \in A_{7}(1243)$ has label $(3,6)$.

The above labeling scheme enables us to derive
 succession rules for $A_{2n+1}(1243)$.
Recall that the functions $f(u)$ and $e(u)$ for a permutation $u \in A_{2n+1}(1243)$  are defined in Theorem \ref{lem:characterize}. In fact, these functions
can be defined on a permutation on any finite set of positive integers. For example,
let $t=48152967$, which is a permutation on $\{1, 2, 4, 5, 6, 7, 8, 9\}$,  we have  $f(t)=7$ and $e(t)=5$.

\begin{thm}\label{lem:1243a}
For $n\geq 0$,
given $u=u_1 u_2 \cdots u_{2n+1}\in A_{2n+1}(1243)$ with label $(a,b)$ ,
the set of labels of the children of  $u$  is given by the set
\[
\{(x,y)\,|\,1 \leq x \leq a+1 , x+2 \leq y \leq b+2\}.
\]
\end{thm}

\pf
  Assume that $w=v_1\mapsto(v_2 \mapsto u)$ is a child of $u$
   and we write $w=w_1w_2\cdots w_{2n+3}$. We aim to
 determine the range of the label $(x,y)$ of $w$.
 According to Theorem \ref{lem:characterize}, we have relations (\ref{rr1}) and (\ref{rr2}), namely,
 $e(u) < v_1 \leq v_2$ and $\text{max}\{u_1+1,f(u)+1\} \leq  v_2 \leq 2n+2$. Since
 $w \in A_{2n+3}(1243)$, from the labeling schemes (\ref{labelscheme1}) and (\ref{labelscheme2}) it follows that
 \begin{align*}
      x&=2n+4-\text{max}\{w_1+1,f(w)+1\},\\[3pt]
      y&=2n+4-e(w).
 \end{align*}
We proceed to compute $f(w)$ and $e(w)$.
Notice that the insertions of $v_1$ and $v_2$ to $u$ may cause new $21$-patterns and $132$-patterns.
Let $s=w_3 w_4 \cdots w_{2n+3}$.
To determine $f(w)$, it suffices to
compare $f(s)$ with the smaller element in each new $21$-pattern.
Similarly, $e(w)$ can be obtained by comparing $e(s)$ with the smallest element in each new $132$-pattern.
Here are two cases.

  \noindent Case 1:  $e(u)<v_1 < v_2$. It is clear that $s$ is order isomorphic to $u$.  We claim that $e(u)=e(s)$.  We first consider the case $e(u)=0$. In this case, by definition
  we see that $u$ is $132$-avoiding. Thus, $s$ is also $132$-avoiding, namely, $e(s)=0$. We now turn to the case
  $e(u)\neq 0$. In other words,
  $e(u)=u_i$ for some $1 \leq i \leq 2n+1$. Since $s$ is order isomorphic to $u$, we deduce that $e(s)=s_i$.
  Hence $u_i=e(u) < v_1 < v_2$. Since $w=v_1 \mapsto (v_2 \mapsto u)$, we see that $s_i=u_i$. This yields $e(s)=e(u)=u_i$. So the claim is verified.

  To compute $e(w)$, we consider the new $132$-patterns caused by
  the insertions of $v_1$ and $v_2$ into $u$.
   Since $w_2=v_2+1$ and $v_2 > f(u)$,
   by the definition of $f(u)$, we find that $w_2$ cannot appear
    as the smallest entry of any $132$-pattern of $w$.
    Hence we need only to consider the new $132$-patterns caused by
    the insertion of $v_1$ into $u$. Since
    $v_1 (v_2+1) v_2$ is a $132$-pattern
    and $v_1>e(u)=e(s)$, we conclude that $e(w)=\text{max}\{v_1, e(s)\}=v_1$.

  To compute $f(w)$, we first determine $f(s)$. There are two cases. If $e(u) < v_1 \leq f(u)$, then $f(u) \neq 0$. Hence $f(u)=u_i$ for some $1 \leq i \leq 2n+1$.
            Since $s$ is order isomorphic to $u$, we have $f(s)=s_i$.
            Noting that $v_2>f(u) \geq v_1$, we get $s_i=u_i+1$. It follows that $f(s)=f(u)+1$. If $f(u) < v_1 < v_2$,
            by the above argument for determining $e(w)$, we have $f(s)= f(u)$. Therefore, both cases we  have $ f(s) \leq f(u)+1$.

  We now consider the new $21$-patterns caused by the insertions of $v_1$ and $v_2$ into $u$.
  Since $v_1<v_2$ and $w_2=v_2+1$, we see that
  $(v_2+1) v_2$ is a $21$-pattern of $w$. Moreover, it can be seen that
  $v_2$ is the largest entry among the smaller elements of the newly formed $21$-patterns. From
  the fact that $ f(s) \leq f(u)+1\leq v_2$ we deduce that
   $f(w)=\text{max}\{v_2, f(s)\}=v_2$.
             Therefore, we have
                \begin{align*}
                 x&=2n+4-\text{max}\{w_1+1,f(w)+1\}=2n+3-v_2,\\[3pt]
                y&=2n+4-e(w)=2n+4-v_1.
                \end{align*}
           Note that $e(u)<v_1< v_2$ and $\text{max}\{u_1+1,f(u)+1\} \leq v_2 \leq 2n+2$,
           we obtain
           \begin{align*}
            &1 \leq x \leq 2n+3-\text{max}\{u_1+1,f(u)+1\},\\[3pt]
            &2n+5-v_2 \leq y \leq 2n+3-e(u).
           \end{align*}
By the labeling rule (\ref{labelscheme1}), namely, $a=2n+2-\text{max}\{u_1+1,f(u)+1\}$, we deduce that $1 \leq x \leq a+1$. By the labeling rule (\ref{labelscheme2}), namely, $b=2n+2-e(u)$ and the fact that $x=2n+3-v_2$,
 we see that $x+2 \leq y \leq b+1$.
 It is easily checked that $(x,y)$ can be any pairs of
 integers such that  $1 \leq x \leq a+1$ and
           $x+2 \leq y \leq b+1$, since
           $e(u)<v_1< v_2$ and $\text{max}\{u_1+1,f(u)+1\} \leq v_2 \leq 2n+2$. This implies that
           the set of labels of the children of $u$ considered in this case
          is given by
           \[\{(x,y)\,|\,1 \leq x \leq a+1~\text{and}~ x+2 \leq y \leq b+1\}.\]

\noindent Case 2: $v_1=v_2$. By (\ref{rr2}),
          we have $v_1=v_2 > f(u)$. Since $s$ is order isomorphic to $u$, by the argument for computing $e(w)$ in Case 1, we find  $f(s)=f(u)$.
           Let us analyze the new $21$-patterns caused by the insertions of $v_1$ and $v_2$ into $u$.
          Clearly, $(v_2+1) (v_2-1)$ is a $21$-pattern of $w$.
          Moreover, it is obvious that $v_2-1$ is the largest entry of the smaller elements in the newly formed  $21$-patterns.
          Since $v_2 -1 \geq f(u)=f(s)$, we have $f(w)=\text{max}\{v_2-1, f(s)\}=v_2-1$.

        By (\ref{rr1}), namely, $e(u) < v_1 \leq v_2$, and  the fact that $s$ is order isomorphic to $u$, it can be seen that $e(s)=e(u)$.
          Since $v_1=v_2 > f(u)$, the insertions of $v_1$ and $v_2$ do not create any new  $132$-patterns.
          It yields that $e(w)=e(s)=e(u)$.
        Consequently, we get
       \begin{align*}
       x&=2n+4-\text{max}\{w_1+1,f(w)+1\}=2n+3-v_2,\\[3pt]
       y&=2n+4-e(w)=2n+4-e(u).
       \end{align*}
       We claim that  $\text{max}\{u_1+1,f(u)+1\} \leq v_1=v_2 \leq 2n+2$. From the definitions of $f(u)$ and $e(u)$,
         it can be seen that $f(u) \geq e(u)$, since each $132$-pattern contains a $21$-pattern. Thus the claim
         follows from (\ref{rr1}) and (\ref{rr2}). This implies that \[ 1 \leq x \leq 2n+3-\text{max}\{u_1+1,f(u)+1\}.\] By the labeling rule (\ref{labelscheme1}), namely, $a= 2n+2- \text{max}\{u_1+1,f(u)+1\}$, we deduce that $1 \leq x \leq a+1$. By the labeling rule (\ref{labelscheme2}), namely, $b=2n+2-e(u)$, we get $y=b+2$.

        Using the same argument as in Case 1, we see that  $(x,y)$  range over all   pairs of
       integers such that  $1 \leq x \leq a+1$ and $y=b+2$.
       Hence the set of labels of the children of $u$ considered in this case is given by
         \[\{(x,y)\,|\,1 \leq x \leq a+1~\text{and}~ y=b+2\}.\]

       Combining Case $1$ and Case $2$,
       the set of labels of the children of  $u$  is given by
        \[\{(x,y)\,|\,1 \leq x \leq a+1, ~x+2 \leq y \leq b+2\},\]
       as required.  This completes the proof.\qed

\begin{figure}
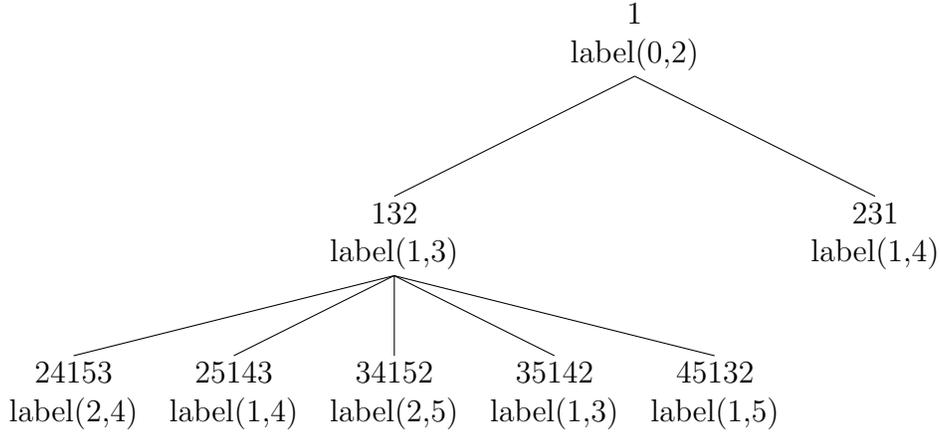

\Tree [.1\\label(0,2) [ 24153\\label(2,4) 25143\\label(1,4) 34152\\label(2,5) 35142\\label(1,3) 45132\\label(1,5) ].132\\label(1,3)  231\\label(1,4) ]
\caption{The first few levels of the generating tree for $A_{2n+1}(1243)$}
 \label{fig:tree1243}
\end{figure}

 Indeed, the above characterization of the labels of the
 children of a permutation $u$ in $A_{2n+1}(1243)$ implies
 that the label of a child $w$ of $u$ is uniquely determined by $w$.
 Thus, in the representation of the generating tree for $ A_{2n+1}(1243)$ we may only keep the labels and ignore
 the alternating permutations themselves. The generating tree
 can be described as follows:
 \begin{align} \label{tree:odd1243}
    \left\{
      \begin{array}{ll}
        root \colon & \hbox{$(0,2)$}, \\ [3pt]
        rule \colon & \hbox{$(a,b)\mapsto \{(x,y)\,| \,1\leq x \leq a+1 ~\text{and}~ x+2 \leq y \leq b+2\}$}.
      \end{array}
    \right.
\end{align}
 Figure \ref{fig:tree1243} gives the first few levels of the generating
 tree for $A_{2n+1}(1243)$.

 Comparing the above description of the generating tree
 for $A_{2n+1}(1243)$ and the generating tree for
 $A_{2n+1}(2143)$ as given by Lewis \cite{Lewis2012}, we arrive at the
 assertion that there is a bijection between    $A_{2n+1}(1243)$ and $A_{2n+1}(2143)$. This proves Theorem \ref{thm:1243=2143}.

The construction of the generating tree
 for $A_{2n+1}(1243)$ can be easily adapted to derive the
generating tree for $A_{2n}(1243)$.
The following theorem gives a similar characterization
 of the set of $1243$-avoiding
alternating permutations in $A_{2n+2}$ that are generated by an
alternating permutation $u$ in $A_{2n}(1243)$.

\begin{thm}\label{lem:characterize1}
For $n \geq 1$, given a permutation $u=u_1u_2\cdots u_{2n} \in A_{2n}(1243)$,
then $w$ is a child of $u$ if and only if
it is of the form $w=v_1 \mapsto (v_2 \mapsto u)$, where
\begin{align*}
&e(u) < v_1 \leq v_2,\\[3pt]
\text{max}\{u_1+1,&
f(u)+1\} \leq  v_2 \leq 2n+1.
\end{align*}
\end{thm}

Based on the above characterization,
we assign a label $(a,b)$ to an alternating
permutation $u=u_1 u_2 \cdots u_{2n} \in A_{2n}(1243)$, where
\begin{align*}
     &a=2n+1-\text{max}\{u_1+1,f(u)+1\}, \\[3pt]
     &b=2n+1-e(u).
\end{align*}
The succession rules for $A_{2n}(1243)$ are exactly the same as these  for $A_{2n+1}(1243)$, since they do not depend on the
parity of the length of $u$.  Note that
$12 \in A_{2n}(1243)$  has label
$(1,3)$. So the generating tree for $ A_{2n}(1243)$ can be described  as follows:
\begin{align}\label{tree:1243}
    \left\{
      \begin{array}{ll}
        root \colon & \hbox{$(1,3)$}, \\[3pt]
        rule \colon & \hbox{$(a,b)\mapsto \{(x,y)\,| \,1\leq x \leq a+1 ~\text{and}~ x+2 \leq y \leq b+2\}$}.
      \end{array}
    \right.
\end{align}

Recall the following generating tree for
$A_{2n}(1234)$ given by Lewis \cite{Lewis2012}:
\begin{align}\label{tree:1234}
    \left\{
      \begin{array}{ll}
        root \colon & \hbox{$(2,3)$}, \\ [3pt]
        rule \colon & \hbox{$(a,b)\mapsto \{(x,y)\,|\,2\leq x \leq a+1 ~\text{and}~ x+1 \leq y \leq b+2\}$}.
      \end{array}
    \right.
\end{align}
Apparently, the above two generating trees are isomorphic via the
 correspondence $(a,b) \rightarrow (a+1,b)$.
This gives another proof of Theorem \ref{thm:1243=1234}, which was proved by B\'ona \cite{BoNA} via
a direct bijection.

\section{Proof of the conjecture $|A_{2n}(4312)|=|A_{2n}(1234)|$} \label{sec:4312}

In this section, we prove the following theorem which leads to the
relation $|A_{2n}(4312)|=|A_{2n}(1234)|$ conjectured by Lewis \cite{Lewis2012}.

\begin{thm}\label{thm:4312=shsyt}
For $n\geq 3$, we have $|A_{2n}(4312)|=|SHSYT(n+2,n,n-2)|$.
\end{thm}

Recall that $|SHSYT(n+2,n,n-2)|$
is known as the $n$-th 3-dimensional Catalan number $C_n^{(3)}$. As
Lewis \cite{Lewis2012} has proved that $|A_{2n}(1234)|=C_n^{(3)}$, we are led to the relation (\ref{c2}),
which we restate as the following theorem.

\begin{thm}\label{thm:4312=1234}
For $n\geq 1$, we have $|A_{2n}(4312)|=|A_{2n}(1234)|$.
\end{thm}

Let us first recall some notation and terminology.
A sequence $\lambda=(\lambda_1,\lambda_2,\cdots,\lambda_m)$ of positive integers is  said to be a partition of $n$ if
 $n=\lambda_1+\lambda_2+\cdots+\lambda_m$ and $\lambda_1 \geq \lambda_2 \geq \cdots \geq \lambda_m >0$, where each $\lambda_i$ is called a part of $\lambda$.
A Young diagram of shape $\lambda$ is defined to be a left-justified array of $n$ boxes with $\lambda_1$ boxes
in the first row, $\lambda_2$ boxes in the second row and so on.
If $\lambda$ is a partition with distinct parts, then the shifted Young
diagram of shape $\lambda$ is an array of cells with $m$ rows, where
each row is indented by one cell to the right with respect to the previous row, and there are $\lambda_i$ cells in row $i$.

A standard Young tableau of shape $\lambda$ is a Young diagram of $\lambda$ whose boxes have been filled with the number $1,2,\ldots, n$
 such that the entries are increasing along each row and each column.
 A shifted standard Young tableau of shape $\lambda$ is a filling of a shifted Young diagram with $1,2,\ldots, n$ such that the
 entries are increasing along each row and each column.
We denote by $SHSYT(\lambda)$ the set of shifted standard Young tableaux of shape $\lambda$.

As shown by Lewis \cite{Lewis2012}, $A_{2n}(1234)$ is
enumerated by the $n$-th $3$-dimensional Catalan number $C_n^{(3)}$ as given in
(\ref{V_n}). To prove relation (\ref{c2}),
it suffices to demonstrate that $A_{2n}(4312)$ is also counted by $C_n^{(3)}$.
In light of the correspondence between $4312$-avoiding alternating permutations and $1243$-avoiding down-up permutations via complement, we proceed to consider the generating tree for $A_{2n+1}(1243)$.

It turns out that the generating tree for $A_{2n+1}(1243)$ is isomorphic to the generating tree for $SHSYT(n+2,n+1,n)$ obtained by Lewis \cite{Lewis2012} as given below:
\begin{align}\label{tree:shsyt}
    \left\{
      \begin{array}{ll}
        root \colon & \hbox{$(1,2)$}, \\ [3pt]
        rule \colon & \hbox{$(a,b)\mapsto \{(x,y)\,|\,2\leq x \leq a+1 ~\text{and}~ x+1 \leq y \leq b+2\}$}.
      \end{array}
    \right.
\end{align}
The above generating tree is based on the following labeling scheme.
Let $T$ be a shifted standard Young tableau $T\in SHSYT(n+2,n+1,n)$,
 and let $T(i,j)$ denote the entry of $T$ in the $i$-th row and the $j$-th column.
 We associate $T$ with a label $(a,b)$,
where $a=3n+4-T(2,n+2)$ and $b=3n+4-T(1,n+2)$.
The isomorphism can be easily established by mapping a label $(a,b)$ in (\ref{tree:shsyt}) to a label $(a-1,b)$ in (\ref{tree:odd1243}).
This leads to the relation
\begin{equation}\label{12shs}
|A_{2n+1}(1243)|=|SHSYT(n+2,n+1,n)|.
\end{equation}

Observe that in the constructions of the generating tree for
$A_{2n+1}(1243)$ and the generating tree for $SHSYT(n+2,n+1,n)$,
the children of any alternating permutation and
any shifted standard Young tableau are uniquely labeled. So the above isomorphism between the generating trees of $A_{2n+1}(1243)$ and $SHSYT(n+2,n+1,n)$ can be restricted to certain classes of
labels. More precisely, we have the following correspondence.

\begin{thm}
For $n\geq 1$,
there is a one-to-one correspondence between
the set $P_n$ of alternating permutations in $A_{2n+1}(1243)$ with labels of the form $(1,b)$ and the
set $Q_n$ of shifted standard Young tableaux in $SHSYT(n+2,n+1,n)$ with labels of the form $(2,b)$.
\end{thm}

By the labeling schemes, we have
\begin{align}
    &P_n=\{u\,|\,u\in A_{2n+1}(1243),~2n+2-\text{max}\{u_1+1,f(u)+1\}=1\},\\[3pt]
    &Q_n=\{T\,|\,T \in SHSYT(n+2,n+1,n),~T(2,n+2)=3n+2\}.
\end{align}
Based on the relation $|P_n|=|Q_n|$, we aim to prove Theorem \ref{thm:4312=shsyt}.
To this end, we shall give characterizations of $P_n$ and $Q_n$  without using the labels.
Then we shall give a bijection between the set $P_n$ and the set $A_{2n}(4312) \cup A_{2n-1}(1243)$
and a bijection between the set $Q_n$ and the set \[ SHSYT(n+1,n,n-1)\cup SHSYT(n+2,n,n-2).\]
In view of (\ref{12shs}), we see that $|A_{2n-1}(1243)|=|SHSYT(n+1,n,n-1)|$. Thus we arrive at Theorem \ref{thm:4312=shsyt}.

\begin{lem} For $n\geq 0$, an alternating permutation $u=u_1u_2\cdots u_{2n+1}\in A_{2n+1}(1243)$
is in $P_n$ if and only if $u_2=2n+1$, that is,
\begin{equation}
P_n=\{u\,|\,u=u_1u_2\cdots u_{2n+1} \in A_{2n+1}(1243),~u_2=2n+1\}.
\end{equation}
\end{lem}

\pf
Assume that  $u=u_1 u_2 \cdots u_{2n+1} \in P_n$, that is, $2n+2-\text{max}\{u_1+1,f(u)+1\}=1$.
It follows that  $u_1=2n$ or $f(u)= 2n$. If $u_1=2n$,
then we have  $u_1 < u_2 \leq 2n+1$ since $u \in A_{2n+1}(1243)$.
This implies that $u_2=2n+1$. If $f(u)= 2n$,  by the definition of $f(u)$,
we find that $2n+1$  precedes  $2n$ in $u$, since $(2n+1)(2n)$ is the only  21-pattern in $u$
with $2n$ being a smaller element.
We aim to show that $u_2=2n+1$. Assume to the contrary that $u_2 < 2n+1$.
Then $u_1 u_2 (2n+1) (2n)$ forms a $1243$-pattern of $u$ since $u_1 < u_2$.
So we have $u_2=2n+1$.

Conversely, assume that $u \in A_{2n+1}(1243)$ and $u_2=2n+1$. We wish to show that
$2n+2-\text{max}\{u_1+1,f(u)+1\}=1$. We consider two cases. If $u_1=2n$, then by the definition of $f(u)$, we have
$f(u)< 2n$. It follows  that $2n+2-\text{max}\{u_1+1,f(u)+1\}=1$. If
$u_1 \neq 2n$, then  $(2n+1)(2n)$ forms a $21$-pattern of $u$ since $u_2=2n+1$.
From the definition of $f(u)$ it can be seen that $f(u)=2n$. Hence we also get
 $2n+2-\text{max}\{u_1+1,f(u)+1\}=1$.
This completes the proof. \qed

\begin{thm}\label{lem:4312a}
For $n\geq 1$, there is a bijection between  $P_n$ and  $A_{2n}(4312) \cup A_{2n-1}(1243)$.
\end{thm}

\pf
We divide $P_n$ into two subsets $P'_n$ and $P''_n$, where
\begin{align*}
P'_n=\{u\,|\, u=u_1 u_2 \cdots u_{2n+1} \in A_{2n+1}(1243),u_2=2n+1 ~\text{and}~u_1> u_3\}, \\[3pt]
P''_n=\{u\,|\, u=u_1 u_2 \cdots u_{2n+1} \in A_{2n+1}(1243),u_2=2n+1 ~\text{and}~ u_1< u_3\}.
\end{align*}
We proceed to show that there is  a bijection
between $P'_n$ and $A_{2n}(4312)$ and there is a bijection between
$P''_n$ and $A_{2n-1}(1243)$.

First, we define a map
$\varphi\colon P'_n \rightarrow A_{2n}(4312)$. Given
 $v=v_1 v_2 \cdots v_{2n+1}\in P'_n$, let $\varphi(v)=\pi^c$,
where $\pi=v_1 v_3v_4 \cdots v_{2n+1}$.
 Note that $\pi$ is a permutation of $[2n]$ since $v_2= 2n+1$.
 Moreover, it can be seen that $\pi$ is $1243$-avoiding since $v$ is $1243$-avoiding and
 $\pi$ is a subsequence of $v$.
From the fact $v_1>v_3<v_4 >\cdots >v_{2n+1}$, we see that
$\pi$ is a down-up permutation. It follows that
 $\pi \in A'_{2n}(1243)$. Hence we deduce that $\varphi(v)=\pi^c \in A_{2n}(4312) $.
That is to say, $\varphi$ is well-defined.

To prove that $\varphi$
is a bijection, we construct the inverse of $\varphi$.
Define a map $\phi\colon A_{2n}(4312)$ $\rightarrow P'_n$. Given $w=w_1w_2\cdots w_{2n}\in A_{2n}(4312)$,
let $\phi(w)=\tau=(2n+1-w_1)(2n+1)(2n+1-w_2)\cdots (2n+1-w_{2n})$.
We claim that $\tau$ is $1243$-avoiding. Since $w$ is $4312$-avoiding,
by complement we see that $\tau_1 \tau_3 \tau_4\cdots \tau_{2n+1}$
is $1243$-avoiding. Note that $\tau_2=2n+1$ does not occur in any $1243$-pattern of $\tau$.
So the claim is verified.  Evidently, $\tau$ is alternating, and hence we have
$\tau \in A_{2n+1}(1243)$.
From the fact that $w_1 <w_2$ we see that $\tau_1 > \tau_3$.
Thus $\tau \in P'_n$, and so $\phi$ is well-defined.
Moreover, it can be easily checked that $\phi= \varphi^{-1}$.
So we conclude that  $\varphi$ is a bijection between  $P'_n$ and $A_{2n}(4312)$.

We next construct a bijection between $P''_n$ and $A_{2n-1}(1243)$.
Given an alternating permutation $u=u_1 u_2 \cdots u_{2n+1}$ in $P''_n$, define
$\psi(u)=st(r)$, where
$r=u_1u_3u_5u_6\cdots u_{2n}u_{2n+1}$ and $st(r)$ is the
permutation of $[2n-1]$ which is order isomorphic to $r$.

We claim that $\psi$ is well-defined, that is,
$\psi(u)$ is an alternating permutation in $A_{2n-1}(1243)$.
Since $u\in P''_n$,  we find that $u_1<u_3<u_4$ and $u_2=2n+1$.
 We assert that $u_3+1=u_4$. Otherwise, $u_1 u_3 u_4 (u_3+1)$ would form a $1243$-pattern of $u$,
contradicting to the fact $u$ is $1243$-avoiding.
Since $u_4>u_5$, we deduce that $u_3>u_5$. It follows that $u_1 < u_3 >u_5 <\cdots <u_{2n}>u_{2n+1}$.
Thus $\psi(u)$ is an alternating permutation of length $2n-1$.
It is clear that
$\psi(u)$ is $1243$-avoiding since $u$ is $1243$-avoiding and $u$ contains $\psi(u)$ as a pattern.
So we conclude that $\psi(u) \in A_{2n-1}(1243)$. This proves the claim.

To prove that $\psi$
is a bijection, we describe the inverse of $\psi$.
Given $q =q_1 q_2 \cdots q_{2n-1}$  in $A_{2n-1}(1243)$,
define
$\theta(q)=p$, where $p=p_1p_2\cdots p_{2n+1}$ is obtained
from $q$ by inserting $2n+1$ after $q_1$ and inserting $q_2+1$ after $q_2$, and increasing each
element of $q$ which is not less than $q_2+1$ by $1$.
For example, for $q=34152 \in A_5(1243)$, we have $p=\theta(q)=3745162$.

We need to show that $\theta$ is well-defined, that is,
$p$ is an alternating permutation in $P''_n$.
By the construction of $p$, we have $p_1=q_1$, $p_2=2n+1$, $p_3=q_2$, $p_4=q_2+1$
and $p_5=q_3$.
It follows that $p_1<p_2>p_3<p_4>p_5$. Since $q_3q_4\cdots q_{2n-1}$ is order isomorphic to
 $p_5p_6\cdots p_{2n+1}$, we find that
$p_5<p_6>\cdots >p_{2n+1}$. This proves that $p$ is alternating.

We proceed to show  that $p$ is $1243$-avoiding. Assume to the contrary
  that $p_tp_ip_jp_k$
forms a $1243$-pattern of $p$.
Since $q$ is $1243$-avoiding, form the construction of $p$, we see
that $p_tp_ip_jp_k$ must contain $p_2$ or $p_4$.
Since $p_2=2n+1$ cannot occur in any $1243$-pattern, $p_4$ appears in $p_tp_ip_jp_k$.
Moreover, $p_3$ must appear in $p_tp_ip_jp_k$. Otherwise, we assume that
$p_3$ does not appear in
$p_tp_ip_jp_k$. Since $p_3=q_2$ and $p_4=q_2+1$, we see that $p_3+1=p_4$.
By replacing $p_4$ with $p_3$ in $p_tp_ip_jp_k$ we obtain a $1243$-pattern
which does not contain $p_4$, a contradiction.
So we have shown that $p_tp_ip_jp_k$ contains both $p_3$ and $p_4$.

Since $p_3+1=p_4$, we have $p_3<p_4$.
By the assumption that $p_tp_ip_jp_k$ forms a $1243$-pattern,
 we have either  $p_tp_i=p_3p_4$ or $p_ip_j=p_3p_4$.
 If $p_tp_i=p_3p_4$, that is, $p_3p_4p_jp_k$ is a $1243$-pattern, where $j>4$.
 Then $p_1p_3p_jp_k$ forms a $1243$-pattern since $p_1<p_3<p_4$, contradicting to the assertion that $p_4$ must appear in any $1243$-pattern of $p$.
 We now consider the case $p_ip_j=p_3p_4$, namely, $p_tp_3p_4p_k$ is a $1243$-pattern, where $k>4$.  This yields that
$p_3<p_k< p_4$. But this is impossible because $p_3+1=p_4$.
This proves that $p$ is $1243$-avoiding.

Till now, we have shown that $p\in A_{2n+1}(1243)$. Combining the fact that
 $p_2=2n+1$ and $p_1<p_3$, we see that $p\in P''_n$.
It follows that $\theta$ is a
well-defined  map from
$A_{2n-1}(1243)$ to $P''_n$.
 It is easy to verify that $\theta=\psi^{-1}$. Hence $\psi$ is a bijection between
$P''_n$ and $A_{2n-1}(1243)$. This completes the proof.
\qed

\begin{thm}\label{lem:4312b}
For $n\geq 3$, there is a bijection between $Q_n$ and
\[SHSYT(n+1,n,n-1)\cup SHSYT(n+2,n,n-2).\]
\end{thm}

\pf
We first decompose $Q_n$ into two subsets $Q_n'$ and $Q_n''$, where
\begin{align*}
   Q'_n&=\{T\,|\,T \in Q_n,~T(3,n+1)=3n+1\},\\[3pt]
   Q''_n&=\{T\,|\,T \in Q_n,~T(1,n+2)=3n+1\}.
\end{align*}
Clearly, $Q'_n \cap Q''_n = \emptyset$.
We wish to show that $Q_n = Q'_n \cup Q''_n $.
It suffices to prove that
$ Q_n \subseteq Q'_n \cup Q''_n$.

Given a shifted standard Young tableau $T$ in $Q_n$, since $T(2,n+2)=3n+2$ and $T(2,n+2)< T(3,n+2) \leq 3n+3$, we find that $T(3,n+2)=3n+3$. Since the entries
in a shifted standard Young tableau
are increasing along each row and each column, we have $T(1,n+2)=3n+1$ or $T(3,n+1)=3n+1$.
So we deduce that $ Q_n \subseteq Q'_n \cup Q''_n$.

We now define a map $\chi$ from
$Q_n$ to the set $SHSYT(n+1,n,n-1)\cup SHSYT(n+2,n,n-2).$
Let $T$ be a shifted standard Young tableau in $Q_n$.
If $T \in Q'_n$,  then let $\chi(T)=T_1$, where $T_1$ is
obtained from $T$ by deleting the boxes $T(2,n+2)$, $T(3,n+1)$ and $T(3,n+2)$.
If $T \in Q''_n$, then let $\chi(T)=T_2$, where $T_2$ is
obtained from $T$ by deleting the boxes $T(1,n+2)$, $T(2,n+2)$ and $T(3,n+2)$.
It is easy to verify that $\chi$ is  well-defined and it is a bijection between
 $Q_n$ and  $SHSYT(n+1,n,n-1)\cup SHSYT(n+2,n,n-2).$
This completes the proof. \qed

 Figure \ref{fig:(2,n+2)=3n+2} gives an illustration
of the two cases when $T(2,n+2)=3n+2$.
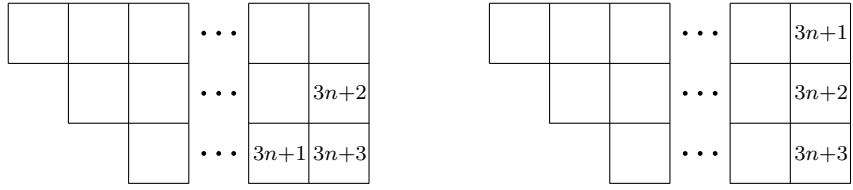
\begin{figure}
\begin{center}\setlength{\unitlength}{.08mm}
\begin{picture}(400,400)(50,0)
\put(-50,0){\line(1,0){200}}
\put(-50,100){\line(1,0){200}}
\put(-50,200){\line(1,0){200}}
\put(-50,300){\line(1,0){200}}
\put(-50,300){\line(0,-1){300}}
\put(50,300){\line(0,-1){300}}
\put(150,300){\line(0,-1){300}}

%\put(0,-100){\circle*{12}} %
\put(-250,0){\line(1,0){100}}
\put(-350,100){\line(1,0){200}}
\put(-450,200){\line(1,0){300}}
\put(-450,300){\line(1,0){300}}
\put(-450,300){\line(0,-1){100}}
\put(-350,300){\line(0,-1){200}}
\put(-250,300){\line(0,-1){300}}
\put(-150,300){\line(0,-1){300}}

\put(-125,250){\circle*{8}}
\put(-100,250){\circle*{8}}
\put(-75,250){\circle*{8}}
\put(-125,150){\circle*{8}}
\put(-100,150){\circle*{8}}
\put(-75,150){\circle*{8}}
\put(-125,50){\circle*{8}}
\put(-100,50){\circle*{8}}
\put(-75,50){\circle*{8}}

\put(857,240){\scriptsize$3n$$+$$1$}
\put(57,140){\scriptsize$3n$$+$$2$}
\put(57,40){\scriptsize$3n$$+$$3$}

\put(750,0){\line(1,0){200}}
\put(750,100){\line(1,0){200}}
\put(750,200){\line(1,0){200}}
\put(750,300){\line(1,0){200}}
\put(750,300){\line(0,-1){300}}
\put(850,300){\line(0,-1){300}}
\put(950,300){\line(0,-1){300}}

%\put(900,-100){\circle*{12}} %
\put(550,0){\line(1,0){100}}
\put(450,100){\line(1,0){200}}
\put(350,200){\line(1,0){300}}
\put(350,300){\line(1,0){300}}
\put(350,300){\line(0,-1){100}}
\put(450,300){\line(0,-1){200}}
\put(550,300){\line(0,-1){300}}
\put(650,300){\line(0,-1){300}}

\put(675,250){\circle*{8}}
\put(700,250){\circle*{8}}
\put(725,250){\circle*{8}}
\put(675,150){\circle*{8}}
\put(700,150){\circle*{8}}
\put(725,150){\circle*{8}}
\put(675,50){\circle*{8}}
\put(700,50){\circle*{8}}
\put(725,50){\circle*{8}}

\put(-43,40){\scriptsize$3n$$+$$1$}
\put(857,140){\scriptsize$3n$$+$$2$}
\put(857,40){\scriptsize$3n$$+$$3$}
\end{picture}
\caption{The two cases when $T(2,n+2)=3n+2$ }
\label{fig:(2,n+2)=3n+2}
\end{center}
\end{figure}

It is clear that
Theorem \ref{thm:4312=shsyt}  follows from Theorem \ref{lem:4312a} and  Theorem \ref{lem:4312b}.
Recall that $|SHSYT(n+2,n,n-2)|$
is counted by the $n$-th 3-dimensional Catalan number $C_n^{(3)}$.
On the other hand, Lewis \cite{Lewis2012} has  shown that $|A_{2n}(1234)|=C_n^{(3)}$.
Thus for $n \geq 3$, we have $|A_{2n}(4312)|=|A_{2n}(1234)|$.
Note that this relation also holds for $n=1,2$.
This completes  the proof of Theorem \ref{thm:4312=1234}.

\section{A generating tree for $A_{2n}(4312)$} \label{sec:4312tree}

In this section, we construct a generating tree for
 $A_{2n}(4312)$. While this generating tree is
  not isomorphic to that for $A_{2n}(1234)$ given by
Lewis \cite{Lewis2012}, it  allows us to give a second proof of
Theorem \ref{thm:1243=2143}, namely, $|A_{2n+1}(1243)|=|A_{2n+1}(2143)|$. To be more specific,
by deleting the leaves of the generating tree for
$A_{2n}(4312)$ and changing every label $(a,b)$ to $(a-1,b)$,
 we are led to the generating  tree for $A_{2n}(3412)$ as given by Lewis \cite{Lewis2012}.
 By restricting this correspondence to certain labels, we obtain Theorem \ref{thm:1243=2143}.

We now present the construction of the generating tree for $A_{2n}(4312)$,
which is analogous to the construction of the generating tree for $A_{2n+1}(1243)$ as given in Section
\ref{sec:1243}. First, we characterize the set of $4312$-avoiding alternating permutations in $A_{2n+2}$ that are generated by an alternating permutation
$u$ in $A_{2n}(4312)$.

\begin{thm} \label{char2n4312}
For $n\geq 1$, given a permutation $u =u_1 u_2 \cdots u_{2n} \in A_{2n}(4312)$, define
\begin{align*}
g(u)=\text{min}\{2n+1,\,u_i\,|\,\text{there exist}~j~\text{and}~k~\text{such that}~i<j<k~ \text{and}~u_j<u_k<u_i\}.
\end{align*}
Then $w$ is a child of $u$ if and only if it is of the form
$w=v_1\mapsto(v_2 \mapsto u)$, where
\begin{equation}\label{v24312}
1\leq v_1\leq v_2.
\end{equation}
and
\begin{equation}\label{v14312}
~u_1+1\leq v_2 \leq g(u),
\end{equation}
\end{thm}

\pf
Suppose that $w=w_1w_2\cdots w_{2n+2}$ is a child of $u$, that is, $w$ is of the form
$v_1\mapsto(v_2\mapsto u)$ and $w\in A_{2n+2}(4312)$. Since $w$ is alternating on $[2n+2]$, we
see that $1\leq v_1\leq v_2\leq 2n+1$ and $v_2\geq u_1+1$. Moreover, we claim that $v_2\leq g(u)$.
Otherwise, there exist $i<j<k$ such that $u_j<u_k<u_i$ and $v_2>u_i$. By the construction of $w$,
we find that $w_2w_{i+2}w_{j+2}w_{k+2}$ forms a $4312$-pattern of $w$, a contradiction.
Hence we are led to the relations  (\ref{v24312}) and (\ref{v14312}).

Conversely, suppose that $v_1$ and $v_2$ are integers satisfying conditions (\ref{v24312}) and (\ref{v14312}).
To prove that $w=v_1\mapsto (v_2 \mapsto u)=w_1w_2\cdots w_{2n+2}$ is a child of $u$,
it suffices  to show that $w$ is an alternating permutation
in $A_{2n+2}(4312)$. Clearly, $w$ is
alternating, since $v_2\geq u_1+1$ and $1\leq v_1\leq v_2$,

It remains to show that $w$ is $4312$-avoiding. Otherwise,  we may
assume that $w_tw_iw_jw_k$ is a $4312$-pattern. We claim that
we can always choose $t=2$. Since $u$ is $4312$-avoiding, by the construction of
$w$, we see that  $w_tw_iw_jw_k$ contains either $w_1$ or $w_2$. If $w_1w_iw_jw_k$
is a $4312$-pattern, since  $w_1 < w_2$ we find that $w_2w_iw_jw_k$ is also a $4312$-pattern.
So the claim is valid.
We continue to prove that  $w_2w_iw_jw_k$ cannot be
 a $4312$-pattern.

Let $s= v_2\mapsto u$
 and write $s=s_1s_2\cdots s_{2n+1}$.
Clearly, $s_1=v_2$ and $s_1s_{i-1}s_{j-1}s_{k-1}$ is a $4312$-pattern of $s$,
since $w_2w_iw_jw_k$ is assumed to be a $4312$-pattern. It follows that $v_2>s_{i-1}=u_{i-2}$.
Note that $u_{i-2}u_{j-2}u_{k-2}$ is a $312$-pattern of $u$. By the definition of $g(u)$,
it can be checked that $u_{i-2}\geq g(u)$. So we get $v_2>u_{i-2}\geq g(u)$,
contradicting to the condition (\ref{v14312}). Hence $w_2w_iw_jw_k$ cannot be
a $4312$-pattern. This implies that $w$ is $4312$-avoiding.
 So we conclude that
 $w$ is an alternating permutation in $A_{2n+2}(4312)$,
 that is to say that $w$ is a child of $u$. This completes the proof. \qed

 Notice  that using the above generating scheme,
 some permutations in $A_{2n}(4312)$ do not have any children.
 Such permutations are called  leaves of the generating tree.
 Permutations having at least one child are called internal vertices
 of the generating tree. For example, the alternating permutation $3412\in A_4(4312)$ is a leaf
 and the alternating permutation $23154867\in A_8(4312)$ is an internal vertex.

 The following two theorems give characterizations
 of leaves and internal vertices of the generating tree for $A_{2n}(4312)$.

\begin{thm}\label{charinterval1}
For $n \geq 1$, given a permutation $u=u_1u_2\cdots u_{2n}\in A_{2n}(4312)$,
$u$ is a leaf
if and only if $g(u)=u_1$.
\end{thm}

\pf We assume that $u$ is a leaf, namely, $u$ has no child.
By Theorem \ref{char2n4312},  we see that if $u$ has a child, then it
 is of the form $v_1\mapsto (v_2\mapsto u)$
satisfying conditions (\ref{v24312}) and (\ref{v14312}), namely, $1\leq v_1\leq v_2$ and $ u_1+1\leq v_2 \leq g(u)$. Now that $u$ has no child, there does not exist  integers $v_1$ and $v_2$
satisfying  (\ref{v24312}) and (\ref{v14312}).
It follows that  $u_1+1 > g(u)$.
Moreover, we claim that  $u_1 \leq g(u)$.
Otherwise, there exist $i <j<k$ such that $u_j <u_k<u_i$ and $u_1 > u_i$.
Consequently, $u_1 u_i u_j u_k$ forms a $4312$-pattern of $u$, which contradicts to the fact that $u$ is
$4312$-avoiding. So the claim is justified. Thus, we have $u_1 \leq g(u) < u_1+1$, namely,
$g(u)=u_1$.

Conversely, assume that $g(u)=u_1$. By Theorem \ref{char2n4312}, it can be easily verified that
$u$ has no child. This completes the proof.
 \qed

\begin{thm}\label{charinterval2}
For $n \geq 1$, given a permutation $u=u_1u_2\cdots u_{2n}\in A_{2n}(4312)$, define
\[h(u)=\text{min}\{2n+1,\,u_j\,|\,\text{there exists }i ~\text{such that}~i<j~\text{and}~u_i < u_j\}.\]
Then $u$ is an internal vertex
if and only if $h(u)=u_1+1$.
\end{thm}

\pf
Assume that $u$ is an internal vertex.
We claim that $u_1 \leq h(u)$. Otherwise, we assume  that $u_1 >h(u)$.
Then there exist $i <j $ such that $u_i <u_j$ and $u_1 >u_j$.
It follows that $u_1 u_i u_j$ forms a $312$-pattern of $u$.
By the definition of $g(u)$, we have $g(u) \leq u_1$.
Meanwhile, since $u$ is $4312$-avoiding, we find that $u_1 \leq g(u)$.
Thus, we reach the equality $u_1= g(u)$. By Theorem \ref{charinterval1},
this  implies that $u$ is a leaf, a contradiction.
Hence  the claim is verified.

Observe that $u_1$ cannot be the second entry of any
$12$-pattern. So by the definition of $h(u)$,  it can be checked that
 $u_1\not=h(u)$. It follows that $u_1 <h(u)$.
On the other hand, since $u_1(u_1+1)$ is a $12$-pattern,  by the definition of $h(u)$,
it can be seen that $h(u)\leq u_1+1$.
In summary, we obtain  $u_1<h(u)\leq u_1+1$, namely, $h(u)=u_1+1$.

Conversely, assume that $h(u)=u_1+1$.
We claim that $h(u) \leq g(u)$.
 If $g(u)=2n+1$, then it is clear that $h(u) \leq g(u)$. If $g(u)< 2n+1$, then there exist $i <j <k$ such that
 $u_j <u_k <u_i$ and $g(u)=u_i$. By the definition of $h(u)$, we see that $h(u) \leq u_k$.
Thus we have $h(u) \leq u_k < u_i = g(u)$.
It follows that for both cases we have $h(u) \leq g(u)$, and so the claim is justified.
 By the assumption that $h(u)=u_1+1$, we obtain $u_1+1\leq g(u)$. By Theorem \ref{char2n4312}, the set of children of $u$ is nonempty. So
 $u$ is an internal vertex.
 This completes the proof.
\qed

To construct the generating tree,  we now give a labeling scheme for
alternating permutations in $A_{2n}(4312)$.
 For $n\geq 1$, given a permutation $u =u_1 u_2 \cdots u_{2n} \in A_{2n}(4312)$,
if $u$ is a leaf, we associate it  with a label $(0,0)$.
If $u$ is an internal vertex,
we associate it with a label $(h(u),g(u))$.
For example, let $u=46253817 \in A_8(4312)$. Since $g(u)=u_1=4$, by Theorem \ref{charinterval1},
we see that $u$ is a leaf. Hence the label of $u$ is $(0,0)$.
It is easily seen that
$12 \in  A_{2}(4312) $ is an internal vertex and it has a label $(2,3)$.

The above labeling scheme enables us to give a
 characterization of the labels of the children generated by $u$.
 Like the extensions of the functions $f(u)$ and $e(u)$ defined in Section \ref{sec:1243} to finite
 integer sequences, the functions $g(u)$ and $h(u)$ can also be extended to
 finite integer  sequences.

\begin{thm} \label{thm1}
Assume that $u=u_1u_2\cdots u_{2n}$ is an alternating permutation in $A_{2n}(4312)$
with label $(a,b)$. If $u$ is an interval vertex, then
it generates ${b-a+1 \choose 2}$ leaves and
the set of labels of the internal vertices generated by $u$ is given by  the set
\begin{align*}
\{(x,y)\,|\,2\leq x\leq a+1,~a+2\leq y\leq b+2\}.
\end{align*}
\end{thm}

\pf
Assume that $w=v_1 \mapsto (v_2 \mapsto u)$ is a child of $u$ and let
$w=w_1w_2\cdots w_{2n+2}$. We aim to characterize the label $(x,y)$ of $w$.
According to Theorem \ref{char2n4312}, we have relations (\ref{v24312}) and (\ref{v14312}), namely, $1 \leq v_1 \leq v_2$
and $u_1+1\leq v_2 \leq g(u)$.
Since $u$ is an internal vertex, it follows from Theorem \ref{charinterval2} that
 $u_1+1=h(u)$. Hence relation (\ref{v14312}) is equivalent to
 $h(u)\leq v_2 \leq g(u)$.

By the labeling scheme, we see that if $w$ is a leaf, then $(x,y)=(0,0)$.
If $w$ is an internal vertex, then $(x,y)=(h(w),g(w))$.
In order to determine the range of $(x,y)$, we distinguish the case when $w$ is a leaf
and the case when $w$ is an internal
vertex. We shall derive expressions of $h(w)$ and $g(w)$ in terms of $v_1, v_2$ and the label $(a,b)$.

Let $s=w_3w_4\cdots w_{2n+2}$. By  the same argument as in the proof of
Theorem \ref{lem:1243a}, we see that in order  to determine $h(w)$, it suffices to compare $h(s)$
with the larger element of each new $12$-pattern caused by the insertions of
$v_1$ and $v_2$ into $u$. The computation of $g(w)$ can be carried out in the same manner.
Here are three cases.

  \noindent Case 1: $h(u)+1\leq v_1\leq v_2, h(u)+1\leq v_2\leq g(u)$.
  By the construction of $w$, we see that $w_1=v_1$ and $w_2=v_2+1$.
  We proceed to judge $w$ is a leaf or an internal vertex. To this end, we
  compute $g(w)$.
  Since both $v_1$ and $v_2$ are not larger than $g(u)$ and $s$ is order isomorphic to $u$, it can be easily verified that  $g(s)=g(u)+2$.

  Now we consider the newly formed $312$-patterns caused by
  the insertions of $v_1$ and $v_2$ into $u$. By the assumption that $h(u)+1\leq v_1\leq v_2$,
  there exist $i <j$ such that $v_1 > u_j >u_i$. It is clear that $w_1w_{i+2}w_{j+2}=v_1u_iu_j$.
  Thus $w_1w_{i+2}w_{j+2}$ forms a $312$-pattern of $w$.
  It is easily verified that $v_1$ is the smallest entry among the
  largest elements of the newly formed $312$-patterns. Since $v_1\leq g(u)$, we deduce that
  $g(w)=\text{min}\{v_1, g(s)\}=v_1=w_1$. By Theorem \ref{charinterval1}, we see that
  $w$ is a leaf. Hence in this case $u$ only generates
  leaves. Using the labeling scheme for $A_{2n}(4312)$, we obtain that $a=h(u)$ and $b=g(u)$.
   So the number of leaves generated by $u$ is given by
\begin{align*}
\sum_{v_2=a+1}^{b}(v_2-a)=1+2+\cdots+(b-a)={b-a+1 \choose 2}.
\end{align*}

  \noindent Case 2: $1\leq v_1\leq h(u), h(u)+1\leq v_2\leq g(u)$.
  Since $s$ is order isomorphic to $u$, using the same argument as in the proof of
   Theorem \ref{lem:1243a}, we obtain that $h(s)=h(u)+1$ and $g(s)=g(u)+2$. To compute $h(w)$,
  we consider the newly formed $12$-patterns caused by the insertions of $v_1$ and $v_2$ into $u$. First,
  $v_1(v_1+1)$ is a newly formed $12$-pattern in $w$. Moreover, it can be seen that $v_1+1$
  is the minimal entry among  the larger elements in the newly formed $12$-patterns.
  Notice that $v_1+1\leq h(u)+1=h(s)$. So we have $h(w)=\text{min}(h(s),v_1+1)=v_1+1$.
  By Theorem \ref{charinterval2}, we find that $w$ is an internal vertex.

  To determine the range of the label of $w$, it suffices  to compute $g(w)$.
  Let us consider the newly formed $312$-patterns in $w$.
  Since $w_1=v_1\leq h(u)$,  we see that $w_1$ cannot occur in any $312$-pattern of $w$.  Since  $v_2\geq h(u)+1$, we deduce that $w_2=v_2+1$ is
  the largest entry of a $312$-pattern in $w$. From the fact that $v_2+1<g(u)+2=g(s)$
  we obtain that $g(w)=\text{min}(v_2+1,g(s))=v_2+1$. Therefore, the label of $w$
  is given by $(x,y)=(v_1+1,v_2+1)$. By the labeling scheme, we see that
  $a=h(u)$ and $b=g(u)$. From the assumption of this case, we get
  $2\leq x\leq a+1$ and $a+2\leq y\leq b+1$.
 This implies that the set of labels of the children of $u$
 considered in this case is given by
\begin{align*}
    \{(x,y)\,|\,2\leq x\leq a+1,~a+2\leq y\leq b+1\}.
\end{align*}

  \noindent Case 3: $1\leq v_1 \leq h(u), v_2=h(u)$.
  Since $s$ is order isomorphic to $u$, using the same argument as
  in the proof of Theorem \ref{lem:1243a},
  we obtain that $h(s)=h(u)+2$. Notice that $w_1(w_1+1)$ is a
  $12$-pattern of $w$ and $w_1+1$ is the minimal entry of the larger elements in the newly formed $12$-patterns caused by the insertions of $v_1$ and $v_2$.
  Since $w_1=v_1\leq h(u)$, we find that
  $h(w)=\text{min}(w_1+1,h(s))=\text{min}(w_1+1,h(u)+2)=w_1+1$. According to Theorem \ref{charinterval2}, $w$ is an internal vertex.

  It remains to determine $g(w)$. Recall that $h(u)\leq g(u)$.
  Hence in this case we have  $v_1 \leq v_2 \leq g(u)$. By the reasoning in the proof of
  Theorem \ref{lem:1243a}, we deduce that $g(s)=g(u)+2$.
  Since $v_1\leq v_2=h(u)$, we see that neither $w_1$ nor $w_2$ can be the largest entry of a $312$-pattern of $w$. This yields that $g(w)=g(s)=g(u)+2$.
  Therefore, the label of $w$ is given by $(x,y)=(v_1+1,g(u)+2)$. Since the label $(a,b)$ is given by
  $a=h(u)$ and $b=g(u)$,
  by the assumption that $1\leq v_1 \leq h(u)$ and  $v_2=h(u)$, we obtain that $2\leq x\leq a+1$ and $y=b+2$.
 It follows that the set of labels of the children of $u$
  considered in this case is given by
\begin{align*}
    \{(x,y)\,|\,2\leq x\leq a+1, y=b+2\}.
\end{align*}

Combining the  above three cases, we see that an internal vertex $u$ generates
${b-a+1 \choose 2}$ leaves and  $a(b-a+1)$ internal vertices labeled by $(x,y)$, where $2\leq x\leq a+1$ and $a+2\leq y\leq b+2$. This completes the proof.
\qed

By Theorem \ref{thm1}, the generating tree for $A_{2n}(4312)$
can be described as follows:
\begin{align*}
    \left\{
      \begin{array}{ll}
        root \colon& \hbox{$(2,3)$},\\ [3pt]
        rule \colon& \hbox{$(a,b)\mapsto \{(x,y)\,|\,2\leq x \leq a+1 ~\text{and}~ a+2 \leq y \leq b+2\} $}\\[8pt]
        &\quad \hbox{$ ~\cup{b-a+1 \choose 2}~ \text{occurrences of}~ (0,0)$.}
      \end{array}
    \right.
\end{align*}

For $n \geq 1$, if we restrict our attention to
 the internal vertices in $A_{2n}(4312)$, then we are led to  the following generating tree:
\begin{align}\label{r4}
    \left\{
      \begin{array}{ll}
        root \colon & \hbox{$(2,3)$}, \\ [3pt]
        rule \colon& \hbox{$(a,b)\mapsto \{(x,y)\,|\, 2\leq x \leq a+1 ~\text{and}~ a+2 \leq y \leq b+2\}$.}
      \end{array}
    \right.
\end{align}
Indeed, the above generating tree is isomorphic to the generating tree for
 $A_{2n}(3412)$ as given by Lewis \cite{Lewis2012}:
\begin{align} \label{r3}
    \left\{
      \begin{array}{ll}
        root \colon & \hbox{$(1,3)$}, \\[3pt]
        rule \colon& \hbox{$(a,b)\mapsto \{(x,y)\,|\, 1\leq x \leq a+1 ~\text{and}~ a+3 \leq y \leq b+2\}$.}
      \end{array}
    \right.
\end{align}

The one-to-one correspondence is easily established by
mapping a label $(a,b)$ in (\ref{r4}) to
a label $(a-1,b)$ in (\ref{r3}). By restricting this correspondence to certain labels,
we arrive at the following bijection.

\begin{thm}\label{one-to-one}
There is a one-to-one correspondence between the set $U_n$ of alternating permutations in $A_{2n}(4312)$ with labels of the form $(2,b)$ and the set
$V_n$ of alternating permutations in $A_{2n}(3412)$ with labels of the form
$(1,b)$.
\end{thm}

The above theorem leads to an alternative proof of Theorem \ref{thm:1243=2143},
that is, for $n \geq 0$, we have $|A_{2n+1}(1243)|=|A_{2n+1}(2143)|$.
To this end, we give characterizations of $U_n$ and $V_n$ without using labels.

\begin{thm}\label{characterizeUnVn}
For $n \geq 1$, we have
\begin{align}
\label{charaU_n}U_n&=\{u\,|\,u=u_1u_2\cdots u_{2n} \in A_{2n}(4312),\, u_1=1\},\\[3pt]
\label{charaV_n}V_n&=\{u\,|\,u= u_1 u_2 \cdots u_{2n}\in A_{2n}(3412),\,u_{2n}=2n\}.
\end{align}
\end{thm}

\pf
Recall that for a permutation $w \in A_{2n}(3412)$  with label $(a,b)$ in the generating tree
defined by Lewis \cite{Lewis2012}, we have $a=d(w)$, where
\begin{equation}
d(w)=2n-\text{max}\{w_i\,|\,\text{there exists }j ~\text{such that}~ j>i ~\text{and}~w_i<w_j\}.
\end{equation}
By Theorem \ref{charinterval2}, a permutation $u=u_1u_2\cdots u_{2n}\in A_{2n}(4312)$ is
an internal vertex if and only if $h(u)=u_1+1$.
Using the labeling schemes for $A_{2n}(4312)$ and $A_{2n}(3412)$, we find that
$U_n$ and $V_n$ can be described in terms of the functions $h(u)$ and $d(u)$, namely,
\begin{align}
U_n&=\{u\,|\,u=u_1u_2\cdots u_{2n}\in A_{2n}(4312),\,h(u)=u_1+1 ~\text{and}~ h(u)=2\},\\[3pt]
V_n&=\{u\,|\,u=u_1u_2\cdots u_{2n}\in A_{2n}(3412),\,d(u)=1\}.
\end{align}

We first prove (\ref{charaU_n}).
Given $u=u_1u_2\cdots u_{2n} \in U_n$, it is easily seen that $u_1=1$.
Conversely, assume that $u=u_1u_2\cdots u_{2n}$ is an alternating permutation
in $A_{2n}(4312)$ with $u_1=1$.
Since the subsequence $12$ forms a $12$-pattern of $u$, by the definition of $h(u)$, we obtain that $h(u)=2$.
Thus the relation $h(u)=u_1+1$ holds.
It follows that $u \in U_n$. Hence  (\ref{charaU_n}) is verified.

We now consider (\ref{charaV_n}). Assume that $u=u_1u_2\cdots u_{2n}$ is an alternating permutation in $V_n$.
Since $d(u)=1$, we see that
$\text{max}\{u_i\,|\,\text{there exists }j ~\text{such that}~ j>i ~\text{and}~u_i<u_j\}= 2n-1$.
Notice that $(2n-1)(2n)$ is the only $12$-pattern of $u$ with $2n-1$ being a smaller element. It follows that $2n-1$ precedes $2n$ in $u$.
If $u_{2n} \neq 2n$, then $(2n-1)(2n)u_{2n-1}u_{2n}$ forms a $3412$-pattern of $u$,
 which is a contradiction. Thus we have $u_{2n}= 2n$. Conversely, if $u_{2n}=2n$, it is easily seen that $d(u)=1$. This completes the
proof. \qed

In view of Theorem \ref{one-to-one}, we see that $|U_n|=|V_n|$.
To prove Theorem \ref{thm:1243=2143}, we shall give a bijection between $U_n$ and $A_{2n-1}(1243)$ and a bijection between  $V_n$ and  $A_{2n-1}(3412)$. Hence Theorem \ref{thm:1243=2143} follows from the
 fact $|A_{2n-1}(2143)|=|A_{2n-1}(3412)|$.

Define a map $\rho \colon U_n \rightarrow A_{2n-1}(1243)$ as follows.
Given an alternating permutation $w= w_1w_2 \cdots w_{2n} \in U_n$,
let $\rho(w)=\pi^c$, where $\pi=(w_2-1) (w_3-1) \cdots (w_{2n}-1)$.
Obviously, $\pi$ is a $4312$-avoiding down-up permutation on $[2n-1]$.
 It follows that $\rho(w) \in A_{2n-1}(1243)$. Thus, $\rho$ is well-defined.
Using the same argument as in the proof of Lemma \ref{lem:4312a},
it can be shown that $\rho$ is a bijection.

To define a map $\mu \colon V_n \rightarrow A_{2n-1}(3412)$,
as we assume that  $u$ is an alternating permutation in $V_n$.
Let $\mu(u)$ be the alternating permutation obtained from $u$
by deleting the last element. It is easy to verify that
$\mu$ is a  bijection.  This gives an alternative proof of Theorem \ref{thm:1243=2143}.

\vspace{0.5cm}
 \noindent{\bf Acknowledgments.}  This work was supported by  the 973
Project, the PCSIRT Project of the Ministry of Education,  and the National Science
Foundation of China.

\end{document}